\documentclass[12pt]{article}

\usepackage{geometry}
\emergencystretch=\maxdimen
\hfuzz=\maxdimen
\usepackage{fourier}  % math font
\usepackage{charter}  % charter
\usepackage{bm}
\usepackage{dsfont}
\usepackage[mathscr]{eucal}

\usepackage{amsmath}
\usepackage{amsfonts}
\usepackage{amsthm}
\usepackage{amssymb}
\usepackage{setspace}
\usepackage{pdflscape}
\usepackage{actuarialsymbol}
\usepackage{tikz}
\usepackage{algorithm}
\usepackage{algorithmicx}

\usepackage{enumitem}
\setenumerate[1]{
	leftmargin=20pt,
	itemsep=-2pt,
	labelsep=0.5em,
	topsep=2pt}
\setitemize[1]{
	leftmargin=15pt,
	itemsep=1pt,
	labelsep=0.5em,
	topsep=2pt}

\usepackage{titling}
\pretitle{\begin{center}\begin{spacing}{1.01}\Large\bfseries}
		\posttitle{\par\end{spacing}
	\end{center}\vspace*{-1.25em}}

\usepackage{titletoc}
\usepackage[explicit]{titlesec}
\titleformat{\section}{\large\bfseries\boldmath\linespread{1.25}\selectfont}{\thesection.}{0.75em}{#1}
\titleformat{\subsection}{\bfseries\boldmath\linespread{0.25}\selectfont}{\thesubsection.}{0.75em}{#1}
\titleformat{\appendix}{\large\bfseries\linespread{1.25}\selectfont}{Appendix\theappendix.}{0.75em}{#1}

\newtheoremstyle{mystyle}{20pt}{20pt}{\itshape}{0cm}{\bfseries}{.}{0.75em}{}
\theoremstyle{mystyle}

\newtheoremstyle{mystyle02}{20pt}{20pt}{}{0cm}{\bfseries}{.}{0.75em}{}
\theoremstyle{mystyle02}
\newtheorem{remark}{Remark}

\makeatletter \@addtoreset{equation}{section} \makeatother
\renewcommand{\theequation}{\arabic{section}.\arabic{equation}}
\usepackage{siunitx}

\usepackage{apacite}
\usepackage{natbib}

\usepackage{lastpage}
\usepackage{fancyhdr}
\usepackage{graphicx}
\usepackage{epstopdf}
\usepackage{subfigure}
\usepackage{overpic}

\usepackage{caption}
\usepackage{color}
\usepackage{hyperref}
\hypersetup{colorlinks=true,linkcolor=blue,citecolor=blue}
\usepackage{bookmark}

\usepackage[hang]{footmisc}
\renewcommand\footnoterule{\kern20pt{\hrule width0.75in height0.5pt}\kern13pt}

\allowdisplaybreaks[4]

\usepackage{caption} %添加图、表的标题
\usepackage{subfigure}
\usepackage{booktabs}
\usepackage{longtable}
\usepackage{multirow}

\begin{document}
	
	\title{Computing the Gerber-Shiu function with interest and a constant 
	dividend barrier by physics-informed neural 
	networks}
	\author{Zan Yu, ~\quad~Lianzeng Zhang\thanks{Corresponding author.}\\[1em]
		{\small School of Finance, Nankai University, Tianjin 300350,  
		China}}
	\date{\vspace*{-3em}}
	\maketitle

	\def\thefootnote{}
	\footnote{E-mail addresses: zhlz@nankai.edu.cn(L. 
		Zhang), yz3006@163.com(Z. Yu)}

	\begin{center}\begin{spacing}{1.05}
			\begin{minipage}[t]{15cm}
				\abstract{
					In this paper, we propose a new efficient method for 
					calculating the Gerber-Shiu discounted penalty function. 
					Generally, the Gerber-Shiu 
					function usually satisfies a class of integro-differential 
					equation. We introduce the physics-informed neural 
					networks (PINN) which embed a differential 
					equation into the loss of the neural network using 
					automatic differentiation. In addition, PINN is more free 
					to set boundary conditions and does not rely on the 
					determination of the initial value. This gives us an idea 
					to calculate more general Gerber-Shiu functions.
					Numerical examples are provided to illustrate the very 
					good performance of our approximation.}\\[7pt]
				\textbf{Keywords} Gerber–-Shiu 
				function, \; Neural networks, \; Integro-differential 
				equation, \; Dividend barrier%\\[7pt]
				%\textbf{Mathematics Subject Classification}~~62P05, 91B30
			\end{minipage}
	\end{spacing}\end{center}
	
	\enlargethispage{-1.75em}
	\thispagestyle{empty}

\section{Introduction} 
In the classical compound Poisson risk model, the surplus process, namely 
$\{U(t)\}_{t \geq 0}$ , has the 
following form:
\begin{equation}
U(t)=u+c t-\sum_{i=1}^{N(t)} X_{i}, \quad t \geq 0,
\end{equation}
where $u \geq 0$ is the initial reserve and $N(t)$ is the number of claims up 
to time $t$ which follows a homogeneous Poisson process of parameter $\lambda 
>0$. The claim sizes  $\left\{X_{i},i=1,2,\ldots\right\}$ are positive, 
independent and identically distributed random variables, with the common 
distribution function $F(x)$ and finite mean $\mu$. We assume that $X_i$ and 
$\{N(t)\}_{t > 0}$ are independent.

If the surplus earns interest at a constant rate $r \geq 0$, the modified 
surplus process can be described by 
\begin{equation}
	U(t)=u \mathrm{e}^{r t}+c \sx*{\angl{t}}[(r)] -\int_0^t 
	\mathrm{e}^{r(t-x)} \mathrm{d} S(x).
	\label{2.2}
\end{equation}
where $S(t)=\sum_{i=1}^{N(t)} X_{i}$ is the aggregate claim amount up to time 
$t$.

The insurers will pay dividends to its shareholders. 
\cite{de1957impostazione} first proposed the dividend barrier strategy for a 
binomial risk model, and since then dividend problems have been widely studied 
under various risk models 
(\cite{lin2003classical}; \cite{albrecher2005distribution}; 
\cite{li2006distribution}; \cite{lin2006compound}; \cite{yuen2007gerber}; 
\cite{gao2008perturbed};
\cite{cai2009analysis}; \cite{gao2010perturbed}; \cite{chi2011threshold};
\cite{zhang2017compound}; \cite{cheung2019periodic} and their references 
therein).

In this paper, we also consider a risk model in which the surplus can earn 
constant interest and dividends are paid according to a constant barrier 
strategy. Within the framework of this barrier strategy, when the surplus 
reaches a barrier of 
constant level $b$, premium income no longer goes into the surplus but is paid 
out as dividends to shareholders. In other words, when the value of the 
surplus reaches $b$, dividends are continuously paid at rate $c+rb$ and the 
surplus remains at level $b$ until the next claim occurs. The incorporation of 
the barrier strategy into \eqref{2.2} yields the capped surplus, 
$\{U_{b}(t)\}_{t \geq 0}$ which can be mathematically expressed by
\begin{equation}
	\mathrm{d} U_b(t)= 
	\begin{cases}
		c \mathrm{d} t-\mathrm{d} S(t)+r U_b\left(t^{-}\right) \mathrm{d} t, 
		& U_b(t)<b, \\ 
		-\mathrm{d} S(t), & U_b(t)=b .
	\end{cases}
	\label{2.3}
\end{equation}

The time of ruin is the first time that the surplus process \eqref{2.3} takes 
a negative value and is denoted by 
$$
T_{b}=\inf \{t \geq 0\mid U_{b}(t)<0\}.
$$
Note that $P\left(T_b<\infty\right)=1$ (i.e., ruin is certain) for $b<\infty$. 
\cite{Gerber1998} first proposed an expected discounted penalty function to 
study 
the time to ruin, the surplus immediately before ruin, and the deficit at ruin 
in the surplus process. Then, the Gerber-Shiu discounted penalty 
function under the surplus process \eqref{2.3} is defined by
\begin{equation}
	\Phi_{b}(u)=\mathbb{E}\left[e^{-\alpha T_{b}} 
	w\left(U_{b}(T_{b}-),\left|U_{b}(T_{b})\right|\right) 
	\mathbf{1}(T_{b}<\infty) \mid U_{b}(0)=u\right], \quad u \geq 0, \alpha 
	\geq 0,
\end{equation}
where $\mathbf{1}(A)$ is the indicator function of event $A$, and 
$w:[0,\infty) \times[0, \infty) \longmapsto[0, \infty)$ is a measurable 
penalty function. 

The Gerber-Shiu function has been widely used by actuarial researchers due to 
its broad applicability in representing various ruin-related quantities. Over 
the past two decades, several stochastic processes have been employed to model 
the temporal evolution of surplus process using the Gerber-Shiu function, such 
as the Sparre-Andersen model 
(\citealt{LANDRIAULT2008600}, \citealt{CHEUNG2010117}), the L{\'e}vy risk 
model 
(\citealt{doi:10.1080/10920277.2006.10597421}, 
\citealt{doi:10.1080/03461238.2011.627747}), and the Markov additive processes 
(\citealt{li_lu_2008}, \citealt{doi:10.1080/10920277.2010.10597599}).
While existing literature has primarily focused on
exploring explicit solutions for Gerber-Shiu functions, this approach has 
limitations as it heavily relies on assumptions about the underlying claim 
size distribution. Under some specific claim size distribution assumptions, 
such as exponential, Erlang and phase-type, explicit formulas are available.
However, for some general individual claim size distributions, it is still 
difficult to obtain explicit expressions. To overcome this challenge, many 
researchers have 
developed some numerical methods to compute the Gerber-Shiu function. For 
example,
\cite{mnatsakanov2008nonparametric} derived the Laplace transform of the 
survival probability, which can then be obtained through numerical inversion 
methods. This method was also applied by \cite{shimizu2011estimation, 
shimizu2012non} to compute the Gerber-Shiu discounted penalty function in
the Lévy risk model and the perturbed model. Other methods have been proposed 
to approximate the Gerber-Shiu function, including truncating Fourier series 
(\citealt{CHAU2015170}, \citealt{zhang_2017, zhang2017estimating}) and using 
the Laguerre basis (\citealt{zhang2018new, zhang2019estimating}).
\cite{wang2019computing} adopted the frame duality projection method to 
compute the Gerber-Shiu function.

Recently, machine learning (ML) has been increasingly used in insurance and 
actuarial science to simulate various financial outcomes. Examples include 
pricing of non-life insurance (\citealt{noll2020case}, \citealt{
wuthrich2023data}), incurred but not reported (IBNR) reserving 
(\citealt{wuthrich2018neural}, \citealt{kuo2019deeptriangle}), the individual 
claims 
simulation (\citealt{gabrielli2018individual}) and mortality forecasting 
(\citealt{hainaut2018neural}).   

This study makes a contribution to the expanding field of research by 
employing physics-informed neural networks (PINNs), a deep-learning technique 
for solving differential equations that has recently been developed and proven 
effective in various studies (\citealt{raissi2019physics}, 
\citealt{lu2021deepxde}, \citealt{karniadakis2021physics}, 
\citealt{yu2022gradient}). Physics-informed neural networks (PINNs) 
incorporate the residuals of differential equations into the loss function of 
the 
neural network using automatic differentiation (\citealt{lu2021deepxde}, 
\citealt{yu2022gradient}). Consequently, the approximation of the solutions 
involves the minimization of the loss function, which can be achieved through  
gradient descent techniques. This mesh-free method of solving equations is 
straightforward to implement and applicable to various types. PINNs have 
achieved success in a diverse range of fields,  including optics 
(\citealt{chen2020physics}), systems biology (\citealt{daneker2023systems}), 
fluid mechanics (\citealt{raissi2020hidden}), solid mechanics 
(\citealt{wu2023effective}). However, there have been fewer applications to 
problems in this study, which this paper will explore.

In this paper, we consider two models. First, we solve the Gerber-Shiu 
function in the absence of a dividend barrier, which we can refer to 
\cite{cai2002expected} for a deep study. Then we consider the discount penalty 
function with a constant dividend barrier, see \cite{yuen2007gerber}. Our 
methodology offers 
significant advantages over traditional methods for solving differential 
equations (DEs). Unlike traditional approaches that involve designing specific 
algorithms for each unique problem, our methodology provides a general 
framework that can be applied to a wide range of DE problems. Compared with 
traditional methods that rely more on initial values, our method is more 
flexible and easier to implement in setting boundary conditions.

The paper is organized as follows. In Section \ref{sec:2}, we present a review 
of 
fundamental findings from the Gerber-Shiu function. In Sections \ref{sec:3} 
and \ref{sec:4}, we 
illustrate the application of the PINN method to solve the Gerber-Shiu 
function. Section \ref{sec:5} provides a summary of our findings and concludes 
the paper.

\section{Preliminaries on risk model and Gerber-Shiu function}\label{sec:2}
Considering $0 \leq u < b$, look at $\tau>0$ such that the surplus has not 
reached level $b$ by time $\tau$, i.e.,$u \mathrm{e}^{r t}+c 
\sx*{\angl{t}}[(r)] 
< b $. By distinguishing according to the time and amount of the first claim, 
if it occurs by $\tau$, and whether the claim causes ruin, we obtain
\begin{equation}
	\begin{aligned}
		\Phi_{b}(u)= & \mathrm{e}^{-(\alpha+\lambda) \tau} \Phi_{b}\left(u 
		\mathrm{e}^{r \tau}+c \sx*{\angl{\tau}}[(r)]\right)+\lambda 
		\int_0^\tau \mathrm{e}^{-(\alpha+\lambda) t} \int_0^{u \mathrm{e}^{r 
		t}+c \sx*{\angl{t}}[(r)]} \Phi_{b} \left(u \mathrm{e}^{r t}+c 
		\sx*{\angl{t}}[(r)]-y\right) \mathrm{d} F(y) \mathrm{d} t \\
		& +\lambda \int_0^\tau \mathrm{e}^{-(\alpha+\lambda) t} A\left(u 
		\mathrm{e}^{r t}+c \sx*{\angl{t}}[(r)]\right) \mathrm{d} t,
	\end{aligned}
    \label{2.4}
\end{equation}
where
$$
A(t)=\int_t^{\infty} w(z, y-z) \mathrm{d} F(y) .
$$
	
By differentiating \eqref{2.4} with respect to $\tau$ and then setting 
$\tau=0$, we obtain the desired integro-differential equation for 
$\Phi_{b}(u)$: 
\begin{equation}
	0=-(\alpha+\lambda) \Phi_{b}(u)+\Phi_{b}^{\prime}(u)(u r+c)+\lambda 
	\int_0^u \Phi_{b}(u-y) \mathrm{d} F(y)+\lambda A(u) .
	\label{2.5}
\end{equation}
	
To obtain a boundary condition, consider $u = b$ and arbitrary $\tau > 0$. 
Analogous to \eqref{2.4}, we have
\begin{equation}
	\Phi_{b}(b)=\mathrm{e}^{-(\alpha+\lambda) \tau} \Phi_{b}(b)+\lambda 
	\int_0^\tau \mathrm{e}^{-(\alpha+\lambda) t} \int_0^b \Phi_{b}(b-y) 
	\mathrm{d} F(y) \mathrm{d} t+\lambda \int_0^\tau 
	\mathrm{e}^{-(\alpha+\lambda) t} A(b) \mathrm{d} t .
	\label{2.6}
\end{equation}
Similarly, by differentiating \eqref{2.6} with respect to $\tau$ and then 
setting $\tau = 0$, we obtain
\begin{equation}
	0=-(\alpha+\lambda)\Phi_{b}(b)+\lambda \int_0^b \Phi_{b}(b-y) \mathrm{d} 
	F(y)+\lambda A(b) .
    \label{2.7}
\end{equation}
Taking limit $u = b$ in \eqref{2.5} and comparing it with \eqref{2.7}, we 
observe that
\begin{equation}
	\Phi_{b}^{\prime}(b)=0
\end{equation}
which implies an intuitive result.

In particular, let $b \rightarrow \infty$,  the discount penalty function 
without constant dividend barrier, $\Phi_{\infty}(u)$, can be expressed by
\begin{equation}
	0=-(\alpha+\lambda) \Phi_{\infty}(u)+\Phi_{\infty}^{\prime}(u)(u 
	r+c)+\lambda 
	\int_0^u \Phi_{\infty}(u-y) \mathrm{d} F(y)+\lambda A(u) .
\end{equation}	
	
\begin{remark}
For \eqref{2.5}, \cite{lin2003classical} and \cite{yuen2007gerber} expressed 
the Gerber-Shiu expected discounted penalty function as the sum of two 
functions. see \ref{appendix:A} for more details. In some special cases with 
exponential claims, we are able to find closed-form expressions for the 
Gerber-Shiu expected discounted penalty function.
\end{remark}
	
\section{Methodology}\label{sec:3}
In this section, we provide a brief overview of deep neural networks, 	
allowing us to present the framework of physics-informed neural networks 
(PINNs), which will be used to solve the integro-differential equation. 
\subsection{Deep neural networks}
In this paper, we consider the classical feed-forward neural networks (FNNs) 
as our deep neural networks. FNNs are easier to train for deep networks.
	
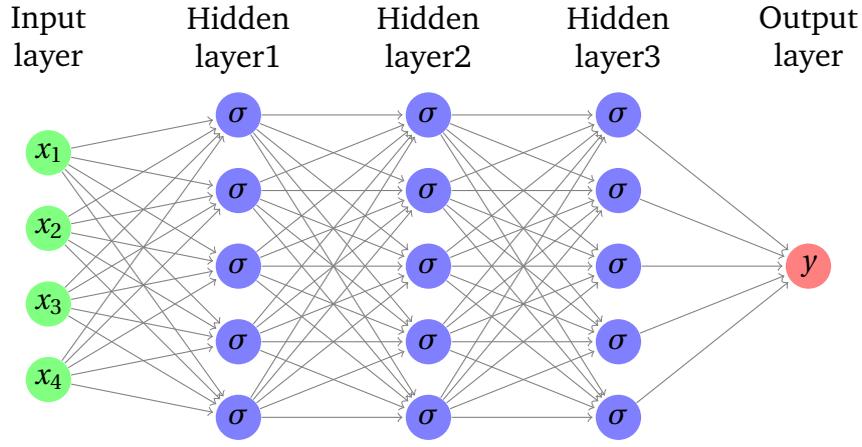
\begin{figure}[h]
	\begin{center}
		\def\layersep{2.5cm}
		\begin{tikzpicture}[shorten >=1pt,->,draw=black!50, node 
		distance=\layersep] 
			\tikzstyle{every pin edge}=[<-,shorten <=1pt]
			\tikzstyle{neuron}=[circle,fill=black!25,minimum size=17pt,inner 
			sep=0pt]
			\tikzstyle{input neuron}=[neuron, fill=green!50];
			
			\tikzstyle{output neuron}=[neuron, fill=red!50];
			
			\tikzstyle{hidden neuron1}=[neuron, fill=blue!50];
			\tikzstyle{hidden neuron2}=[neuron, fill=blue!50];
			\tikzstyle{hidden neuron3}=[neuron, fill=blue!50];
			
			\tikzstyle{annot} = [text width=4em, text centered]
				
% Draw the input layer nodes
			\foreach \name / \y in {1,...,4}
% This is the same as writing \foreach \name / \y in {1/1,2/2,3/3,4/4}
			\node[input neuron] (I-\name) at (0,-\y) {$x_{\y}$};
				
% Draw the hidden layer1 nodes
			\foreach \name / \y in {1,...,5}
			\path[yshift=0.5cm] node[hidden neuron1] (H1-\name) at 
			(\layersep,-\y cm) {$\sigma$};
%Draw the hidden layer2 nodes
			\foreach \name / \y in {1,...,5} 
			\path[yshift=0.5cm] node[hidden neuron2,right of=H1] (H2-\name) at 
			(\layersep,-\y cm){$\sigma$};
%Draw the hidden layer3 nodes
			\foreach \name / \y in {1,...,5}
			\path[yshift=0.5cm] node[hidden neuron3,right of=H2] (H3-\name) at 
			(2*\layersep,-\y cm){$\sigma$};
				
% Draw the output layer node
			\node[output neuron, right of=H3-3] (O) {$y$};
				
% Connect every node in the input layer with every node in the hidden layer.
			\foreach \source in {1,...,4}
			\foreach \dest in {1,...,5}
			\path (I-\source) edge (H1-\dest);

			\foreach \source in {1,...,5}
			\foreach \dest in {1,...,5}
			\path (H1-\source) edge (H2-\dest);
				
			\foreach \source in {1,...,5}
			\foreach \dest in {1,...,5}
			\path (H2-\source) edge (H3-\dest);
% Connect every node in the hidden layer with the output layer
			\foreach \source in {1,...,5}
			\path (H3-\source) edge (O);
				
% Annotate the layers
			\node[annot,above of=H1-1, node distance=1cm] (hl) {Hidden layer1};
			\node[annot,left of=hl] {Input layer};
			\node[annot,right of=hl] (h2){Hidden layer2};
			\node[annot,right of=h2] (h3){Hidden layer3};
			\node[annot,right of=h3] {Output layer};
			\end{tikzpicture}
	\end{center}
\caption{\textbf{Visualization of a deep neural network.} In this example, the 
number of layers $L$ is $4$. In other words, the network consists of an input 
layer, three hidden layers, and an output layer.}
\label{fig:1}
\end{figure}

Let $\mathcal{N}^L(\mathbf{x}): \mathbb{R}^{d_{\text {in}}} \rightarrow 
\mathbb{R}^{d_{\text {out}}}$ be an $L-$layer neural network with $L-1$ hidden 
layers. The $\ell$-th layer has $N_{\ell}$ neurons, where $N_{0}=d_{\text 
{in}}$ and $N_{L}=d_{\text {out}}$.  For each $1 \leq \ell \leq L$, we define 
a weight matrix $\boldsymbol{W}^{\ell} \in \mathbb{R}^{N_{\ell} \times 
N_{\ell-1}}$ and bias vector $\mathbf{b}^{\ell} \in \mathbb{R}^{N_{\ell}}$ in 
the $\ell$th layer, respectively. Given a nonlinear activation function 
$\sigma$, the FNN is recursively defined as follows:
\begin{equation*}
	\begin{aligned}
		\text{input layer:}&  \quad  \mathcal{N}^0(\mathbf{x})=\mathbf{x} \in 
		\mathbb{R}^{d_{\mathrm{in}}}, \\
		\text{hidden layers:}&  \quad 
		\mathcal{N}^{\ell}(\mathbf{x})=\sigma\left(\boldsymbol{W}^{\ell} 
		\mathcal{N}^{\ell-1}(\mathbf{x})+\boldsymbol{b}^{\ell}\right) \in 
		\mathbb{R}^{N_{\ell}} \quad \text{for} \quad 1 \leq \ell \leq L-1, \\
		\text{output layer:}&  \quad 
		\mathcal{N}^L(\mathbf{x})=\boldsymbol{W}^L 
		\mathcal{N}^{L-1}(\mathbf{x})+\boldsymbol{b}^L \in 
		\mathbb{R}^{d_{\text {out}}} ;
	\end{aligned}
\end{equation*}
A network with $L = 4$ is visualized in Figure \ref{fig:1}. There are 
different 
possible activation functions $\sigma$, and in this study, we use the 
hyperbolic tangent (tanh), defined as
$$
\tanh(x)=\frac{e^x-e^{-x}}{e^x+e^{-x}}.
$$

\subsection{Physics-Informed Neural Networks (PINN)} 
In this part, we introduce Physics-Informed Neural Networks (PINN) and 
explain that how it can be used to solve differential equations. We start by 
considering a generic single differential equation:
\begin{equation}
		\begin{aligned}
			& \mathcal{D}(u(\boldsymbol{x}))=f(\boldsymbol{x}), \quad 
			\boldsymbol{x} \in \Omega \\
			& \mathcal{B}(u(\boldsymbol{x}))=g(\boldsymbol{x}), \quad 
			\boldsymbol{x} \in \partial \Omega
		\end{aligned}
\end{equation}
where $\mathcal{D}(\cdot)$ and $\mathcal{B}(\cdot)$ represent some linear or 
nonlinear differential operators, $u: \mathbb{R}^n \rightarrow \mathbb{R}$ is 
the unknown exact solution to the differential equation, 
$\boldsymbol{x}=\left\{x_1, \ldots, x_n\right\}^T$ are the independent 
variables (with $x_i \in \mathbb{R}, \forall i=1, \ldots, n$ ), $f: 
\mathbb{R}^n \rightarrow \mathbb{R}$ and $g: \mathbb{R}^n \rightarrow 
\mathbb{R}$ are some known functions, $\Omega$ is a bounded domain with 
boundary $\partial \Omega$.
	
In the PINN framework, as presented by \cite{RAISSI2019686}, the function $u$ 
is approximated using a neural network 
$\hat{u}(\boldsymbol{x};\boldsymbol{\theta})$, where $\boldsymbol{\theta}$ 
represents the parameters (e.g., weights and biases) of the neural network 
$\hat{u}$ . 
The network $\hat{u}$ takes the input $\boldsymbol{x}$ and outputs a vector 
with the same dimension as $u$. By utilizing automatic differentiation in 
popular machine learning frameworks like Tensorflow or PyTorch, we can easily 
compute the derivatives of $\hat{u}$ with respect to its inputs 
$\boldsymbol{x}$ using the 
chain rule for differentiating compositions of functions.
	
To satisfy the constraints of the differential equation and boundary 
conditions, we impose restrictions on the neural network $\hat{u}$. This is 
achieved by constraining $\hat{u}$ at scattered points, which make up the 
training data. The training data consists of two sets: 
$\left\{\boldsymbol{x}_{\Omega}\right\} \subset \Omega$, representing points 
in the domain, and $\left\{\boldsymbol{x}_{\partial \Omega}\right\} \subset 
\partial\Omega$, representing points on the boundary. These sets are referred 
to as the 'residual points'. To evaluate the disparity between the neural 
network approximation $\hat{u}$ and the given constraints, we define a loss 
function that consists of the weighted sum of the $L^2$ norms of the residuals 
for both the equation and the boundary conditions. Mathematically, this can be 
expressed as follows:
	
\begin{equation}
	\mathcal{L}(\boldsymbol{\theta}) = 
	w_{f}\sum_{x\in\left\{\boldsymbol{x}_{\Omega}\right\}}
	\left(\left\|\mathcal{D}\left(\hat{u}\left(x;\boldsymbol{\theta}\right)\right)
	-f\left(x\right)\right\|^2\right)
	+ w_{g}\sum_{x\in\left\{\boldsymbol{x}_{\partial\Omega}\right\}}
	\left(\left\|\mathcal{B}\left(\hat{u}\left(x;\boldsymbol{\theta}\right)\right)
	-g\left(x\right)\right\|^2\right),
\end{equation}
where $w_{f}$ and $w_{g}$ denote the weights assigned to the equation and 	
boundary conditions, respectively. In order to obtain an optimal solution 
$\boldsymbol{\theta}^*$, we aim to minimize the loss function 
$\mathcal{L}(\boldsymbol{\theta})$ using gradient-based optimization 
techniques, such as Adam optimizer (\cite{kingman2015adam}) or L-BFGS 
optimizer (\cite{byrd1995limited}). Thus, our optimal solution can be 
expressed as:
	
\begin{equation}
	\boldsymbol{\theta}^* = \arg \min _{\boldsymbol{\theta}} 
	\mathcal{L}(\boldsymbol{\theta}).
\end{equation}
	
Finally, we summarize the above procedure in Algorithm \ref{algorithm:1}.	
\begin{algorithm}[htb!]
	Step 1: Initialize the parameters $\boldsymbol{\theta}$ of the neural 
	network $\hat{u}(\boldsymbol{x};\boldsymbol{\theta})$.\\
	Step 2: Set up the training sets 
	$\left\{\boldsymbol{x}_{\Omega}\right\}$ 
	and $\left\{\boldsymbol{x}_{\partial \Omega}\right\}$ for the equation 
	and 
	boundary/initial conditions, respectively.\\
	Step 3: Define the loss function by calculating the weighted $L^2$ 
	norm of 
	both the differential equation and boundary condition residuals.\\
	Step 4: Train the neural network by minimizing the loss function 
	$\mathcal{L}(\boldsymbol{\theta})$ to obtain the optimal parameters 
	$\boldsymbol{\theta}^*$.
	\caption{Physics-Informed Neural Network (PINN) algorithm for solving 
	differential equations.}
	\label{algorithm:1}
\end{algorithm}
	
It is worth mentioning that in this study, when solving integro-differential 
equations (IDEs), the automatic differentiation technique is employed to 
analytically derive the integer-order derivatives. However, when dealing with 
the integral term on the right side of \eqref{2.5}, traditional methods 
such as Gaussian quadrature are used to discretize the integral operators. The 
following discrete form can be represented as follows:
$$
\begin{aligned}
	\int_0^u \Phi_{b}(u-y) \mathrm{d} F(y)&=\int_0^u f(u-y) \Phi_{b}(y) 
	\mathrm{d}y\\
	&\approx \sum_{i=1}^n w_i f({u-y_i}) \Phi_{b}(y_i)
\end{aligned}
$$ 	
Then, a PINN is used to solve the following differential equations instead of 
the original equation:
\begin{equation}
	0=-(\alpha+\lambda) \Phi_{b}(u)+\Phi_{b}^{\prime}(u)(u r+c)+\lambda 
	\sum_{i=1}^n w_i f({u-y_i}) \Phi_{b}(y_i)
	+\lambda A(u) .	
\end{equation}

\section{Numerical illustrations}\label{sec:4}
In this section, we present some numerical examples to show that the PINN is 
very efficient for computing the Gerber-Shiu function. All results are 
performed in Python on Windows, with Intel(R) Core(TM) i7 
CPU, at 2.60 GHz and a RAM of 16 GB. In all examples, we use the tanh as the 
activation function, and the other hyperparameters are listed in Table 
\ref{tab:1}. The 
weights $w_f$ , $w_g$ in the loss function are set to $1$.

\begin{table}[htbp]
	\centering
	\caption{Hyperparameters used for all the following examples. The
		optimizer L-BFGS does not require learning rate, and the neural 
		network (NN) is trained until convergence, therefore the number of 
		iterations 
		is also ignored for L-BFGS.}
	\begin{tabular}{cccc}
		\toprule
		NN depth  & NN width & Optimizer &  \#Iterations \\
		\midrule
	20	&    4   &  L-BFGS     & - \\
		\bottomrule
	\end{tabular}%
	\label{tab:1}%
\end{table}%

Throughout this section, we set $c = 
1.5$, $\lambda = 1$, $r=0.01$ and we consider the following three claim size 
densities:
\begin{enumerate}
	\item[(1)] Exponential density: $f(x)=\mathrm{e}^{-x}, x>0$.
	\item[(2)] Erlang(2) density: $f(x)=4x\mathrm{e}^{-2x}, x>0 $.
	\item[(3)] Combination of exponentials density: $f(x)=3 \mathrm{e}^{-1.5 
		x}-3 \mathrm{e}^{-3 x}$.
\end{enumerate}
Note that $\mathbb{E}(X) = 1$ under the above three claim size density 
assumptions. We 
will compute the following four special Gerber-Shiu functions:  
\begin{enumerate}
	\item[(1)] Ruin probability: $ 
	\Phi(u)=\mathbb{P}\left(T<\infty \mid 
	U(0)=u\right)$, where $\alpha = 0$, $ w(x, y) \equiv 1$.
	\item[(2)] Laplace transform of ruin time 
	$\Phi(u)=\mathbb{E}\left[e^{-\alpha T}\mathbf{1}(T<\infty) \mid 
	U(0)=u\right]$, where $\alpha = 0.01$,  $w(x, y) \equiv 1$.
	\item[(3)] Expected claim size causing ruin: 
	$\Phi(u)=\mathbb{E}\left[\left(U(T-)+\left|U(T)\right|\right) 
	\mathbf{1}(T<\infty) \mid U(0)=u\right]$, where $\alpha = 0$, $w(x, 
	y)=x+y$.
	\item[(4)] Expected deficit at ruin: 
	$\Phi(u)=\mathbb{E}\left[\left|U(T)\right| \mathbf{1}(T<\infty) \mid 
	U(0)=u\right]$, where $\alpha = 0$, $w(x, y)=y$.
\end{enumerate}

\begin{figure}[h]
	\centering  %图片全局居中
	\subfigure[]{
		\label{Fig.sub.1}
		\includegraphics[width=0.45\textwidth]{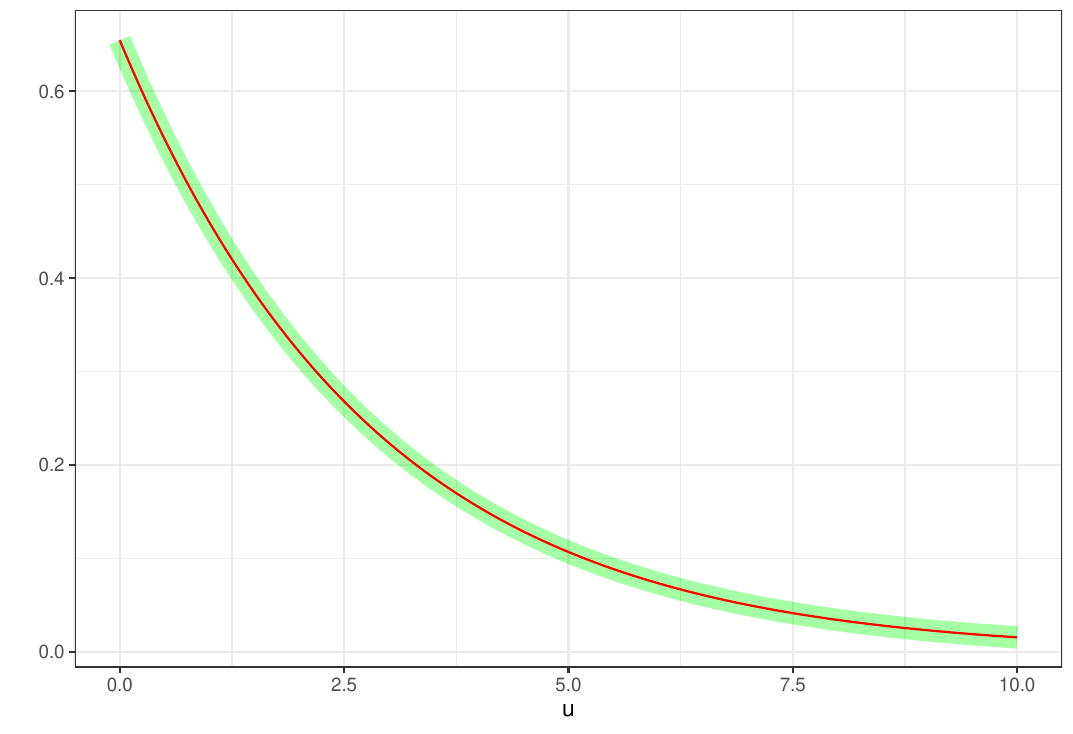}}
	\subfigure[]{
		\label{Fig.sub.2}
		\includegraphics[width=0.45\textwidth]{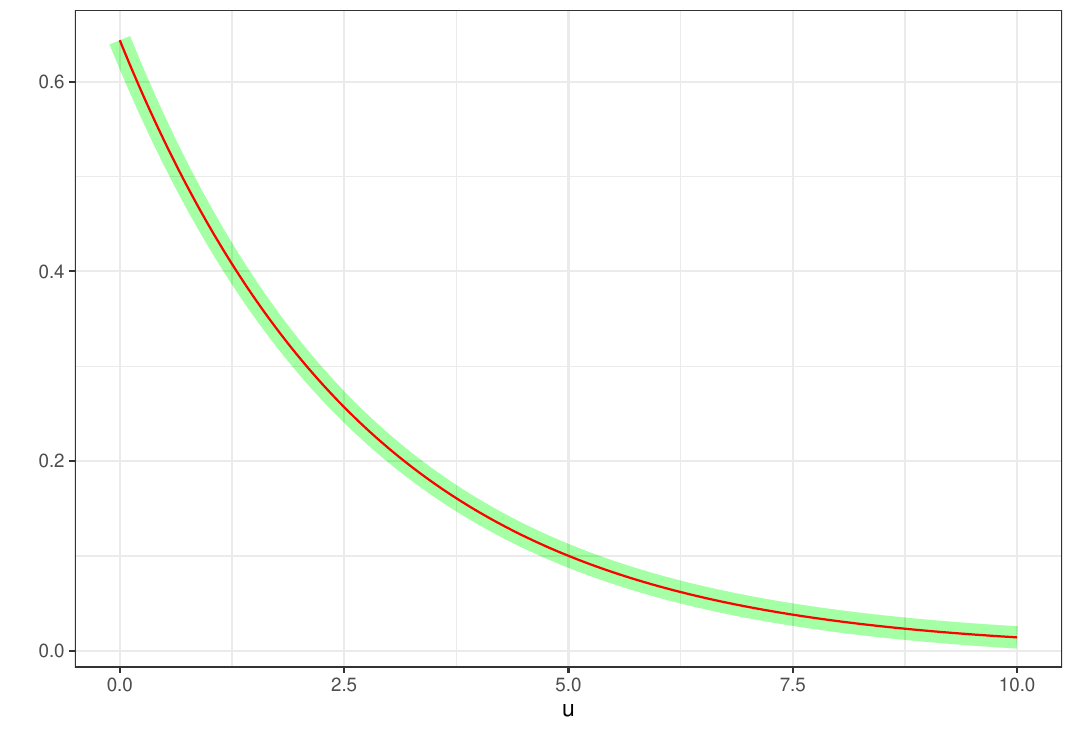}}
	\subfigure[]{
		\label{Fig.sub.3}
		\includegraphics[width=0.45\textwidth]{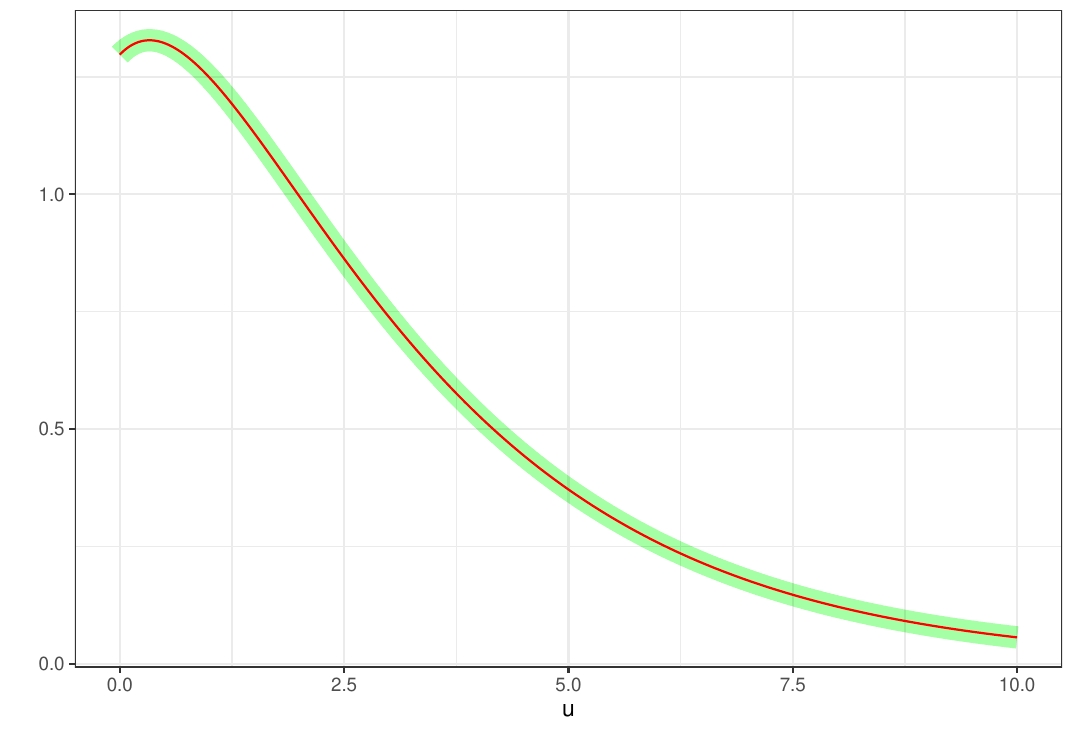}}
		\subfigure[]{
		\label{Fig.sub.4}
		\includegraphics[width=0.45\textwidth]{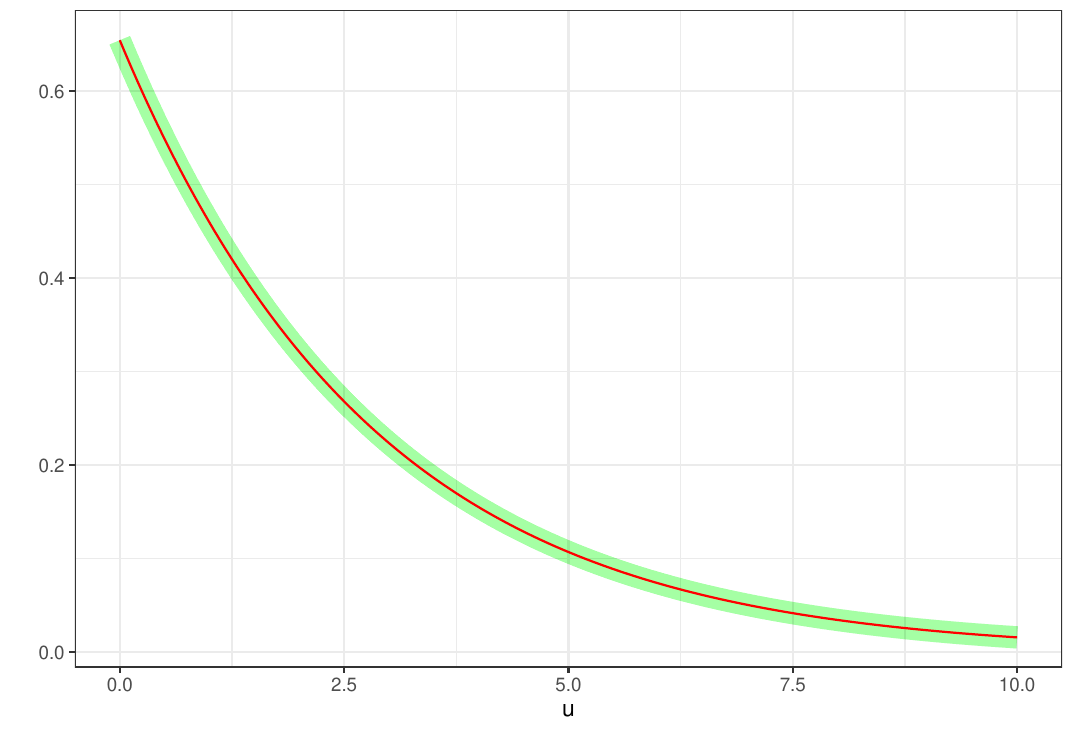}}
	\caption{Computing the Gerber-Shiu functions with Exponential claim size 
		density where red lines
		denote the value curves using collocation method and green lines 
		denote estimate value curves. (a) ruin probability; (b) Laplace 
		transform of ruin time; (c) expected claim size 
		causing ruin; (d) expected deficit at ruin.}
	\label{Fig.2}
\end{figure}

\begin{figure}[h]
	\centering  %图片全局居中
	\subfigure[]{
		\label{Fig.sub.5}
		\includegraphics[width=0.45\textwidth]{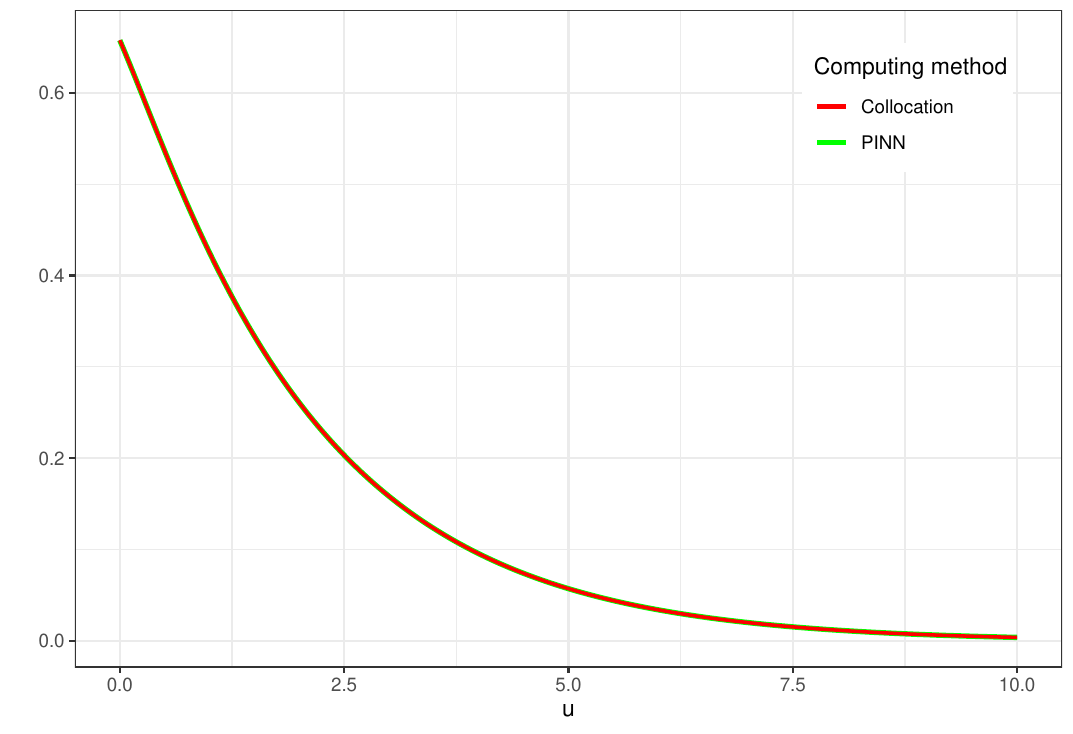}}
	\subfigure[]{
		\label{Fig.sub.6}
		\includegraphics[width=0.45\textwidth]{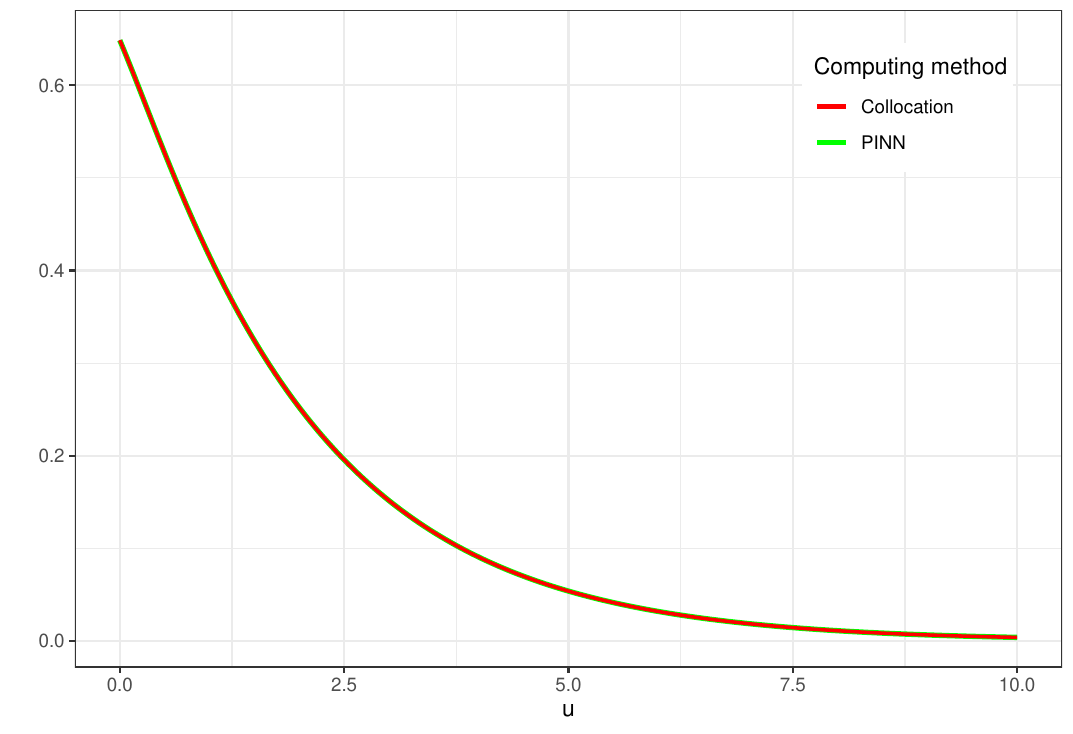}}
	\subfigure[]{
		\label{Fig.sub.7}
		\includegraphics[width=0.45\textwidth]{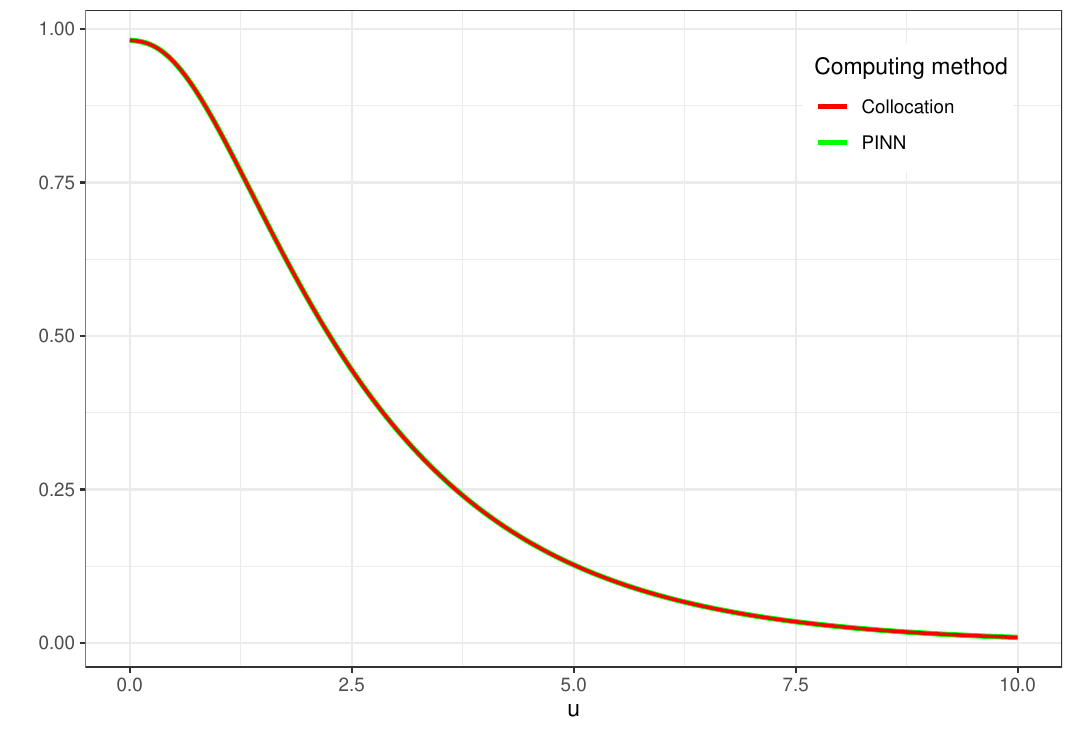}}
	\subfigure[]{
		\label{Fig.sub.8}
		\includegraphics[width=0.45\textwidth]{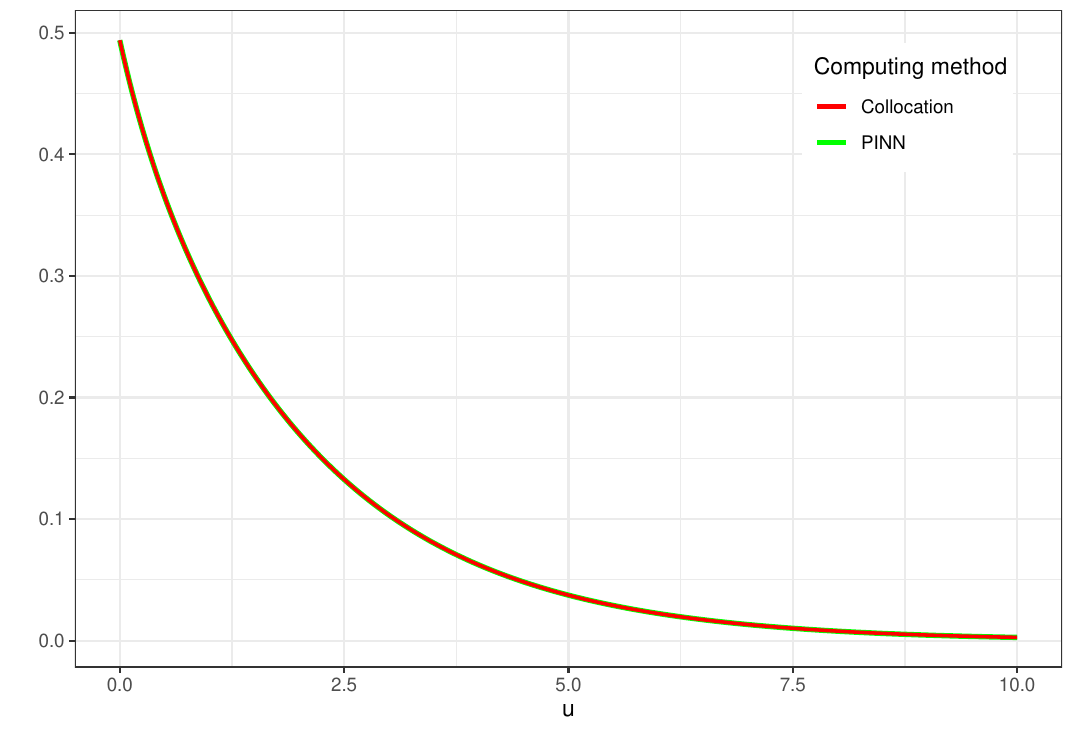}}
	\caption{Computing the Gerber-Shiu functions of Erlang (2) claim sizes 
		density by mean
		value curves. (a) ruin probability; (b) Laplace 
		transform of ruin time; (c) expected claim size 
		causing ruin; (d) expected deficit at ruin.}
	\label{Fig.3}
\end{figure}

\begin{figure}[h]
	\centering  %图片全局居中
	\subfigure[]{
		\label{Fig.sub.9}
		\includegraphics[width=0.45\textwidth]{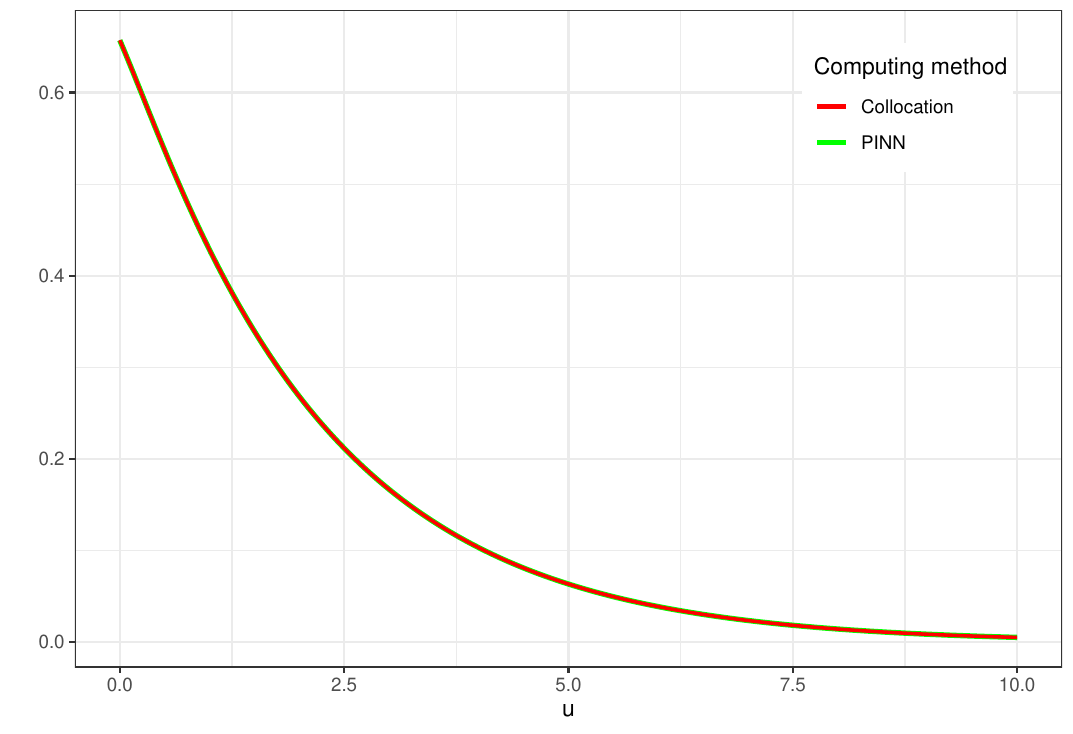}}
	\subfigure[]{
		\label{Fig.sub.10}
		\includegraphics[width=0.45\textwidth]{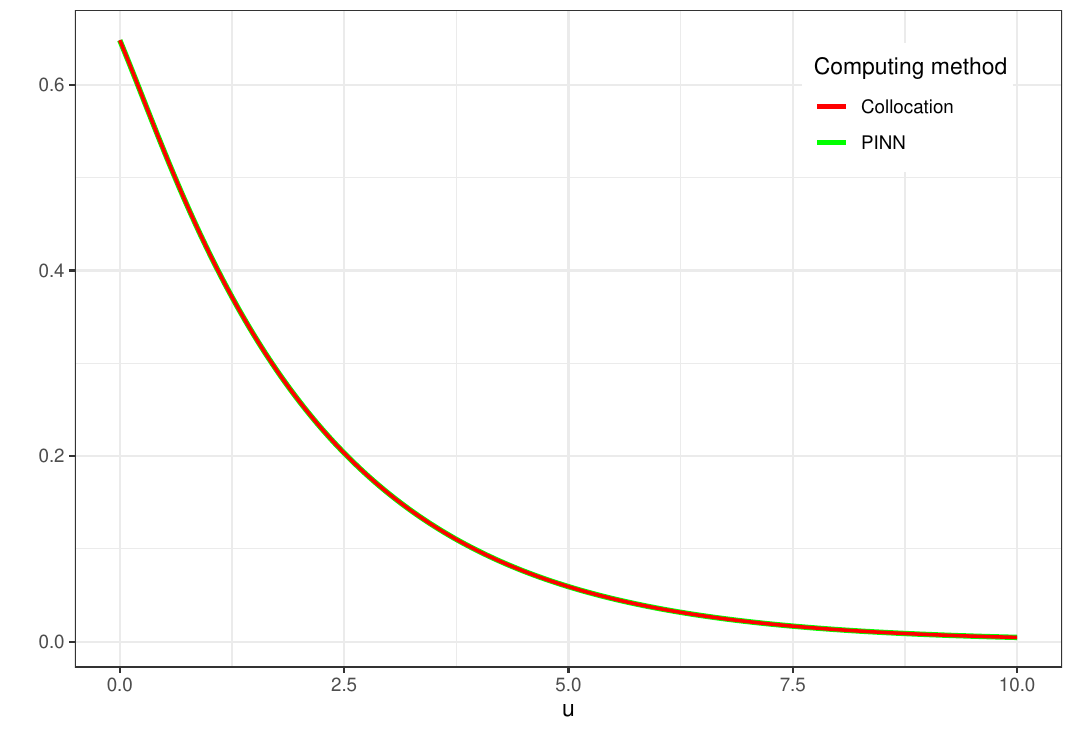}}
	\subfigure[]{
		\label{Fig.sub.11}
		\includegraphics[width=0.45\textwidth]{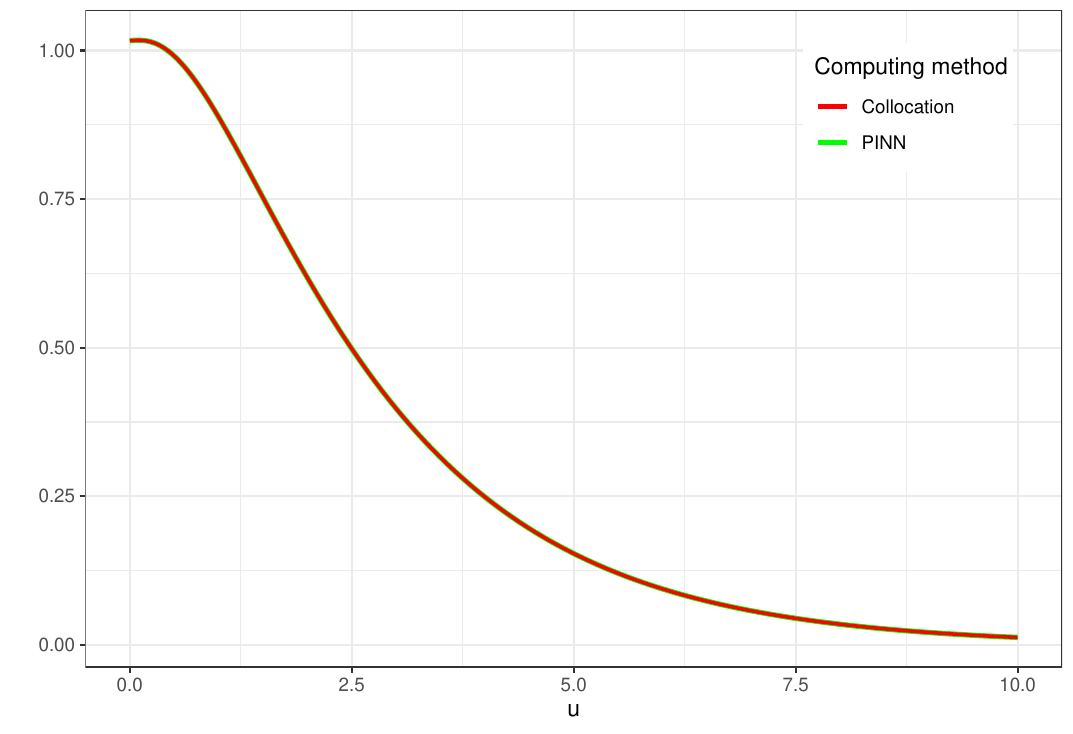}}
	\subfigure[]{
		\label{Fig.sub.12}
		\includegraphics[width=0.45\textwidth]{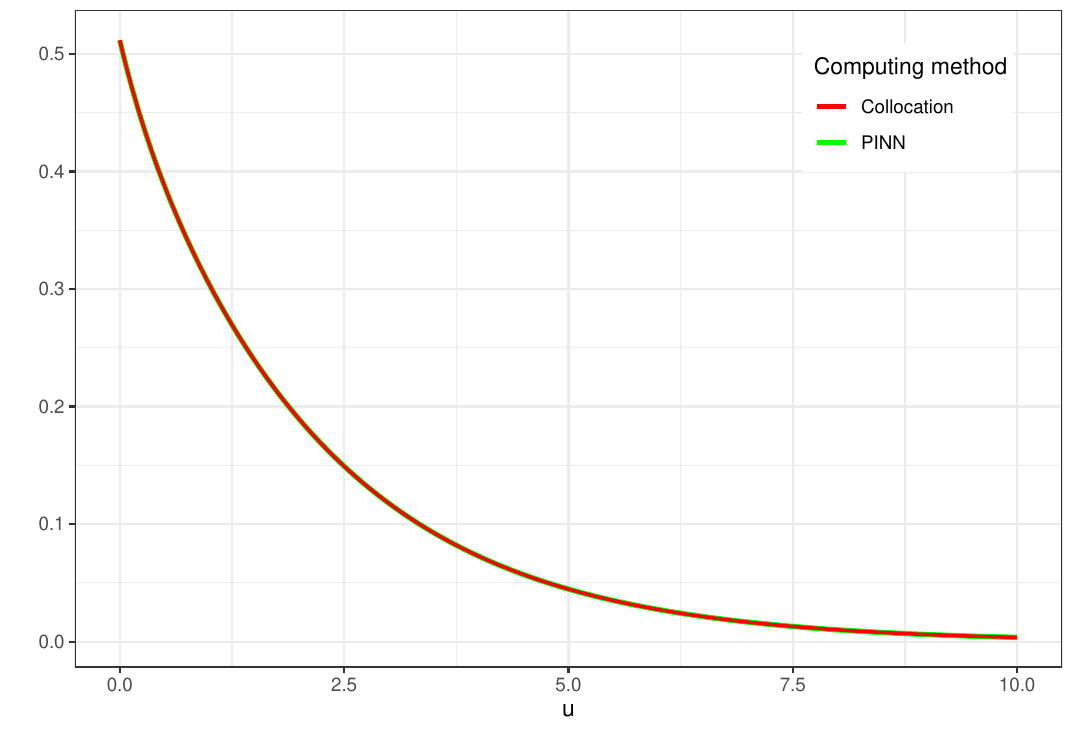}}
	\caption{Computing the Gerber-Shiu functions of 
	combination-of-exponentials claim sizes 
		density by mean
		value curves. (a) ruin probability; (b) Laplace 
		transform of ruin time; (c) expected claim size 
		causing ruin; (d) expected deficit at ruin}
	\label{Fig.4}
\end{figure}

First, we calculate the Gerber-Shiu function without the dividend barrier. For 
more general Gerber-Shiu functions, it is hard to find the
explicit results. So we compare the results of our calculations with classical 
numerical methods. $\Phi_{\infty}(0)$ is determined by \ref{appendix:B}.
Since the results of neural networks closely depend on the selection of random 
seeds, we plot the estimates from 20 experiments on the same figure to show 
variability bands and illustrate the stability of the procedures in Figure 
\ref{Fig.2}. We observe 
that the results are very close to each other. Moreover, they are close to the 
values calculated by collocation method (\citealt{yu2023computing}). Next, for 
Erlang (2) claim size 
density and combination of exponentials claim size density,
we compare the mean value curves of the estimated Gerber-Shiu functions with  
the values computing by collocation method with $N=512$ in Figures 
\ref{Fig.3}--\ref{Fig.4}. Obviously, the mean value curves perform very well, 
and it is difficult to distinguish them from each other. Then we compute the  
mean relative error for all three distribution functions. Table \ref{tab:2} 
shows that the maximum relative error between these two calculation methods. 
Note that the maximum relative error is about \num{1e-5}, which also shows 
that our method performs well.
  
\begin{table}[h]
	\centering
	\caption{Maximum relative error}
		\resizebox{\textwidth}{!}{
	\begin{tabular}{lcccc}
		\toprule
		Claim size & Ruin probability & Laplace transform of ruin time & 
		Expected claim size causing ruin & expected deficit at ruin \\
		\midrule
		Exponential & 3.1201E-05 & 2.3182E-05 & 3.9432E-05 & 1.5646E-05 \\
		Erlang (2) & 3.8662E-05 & 4.6637E-05 & 4.4951E-05 & 5.5839E-05 \\
		Combination of exponentials & 4.1481E-05 & 4.9122E-05 & 8.0658E-05 & 
		2.6430E-05 \\
		\bottomrule
	\end{tabular}}
	\label{tab:2}%
\end{table}% 

Then we compute the Gerber-Shiu function with a constant dividend barrier, 
here the boundary condition of $\Phi_{b}^{\prime}(b)=0$ is only set without 
the need to calculating the initial value, $\Phi_{b}(0)$. Here we also set $c 
= 1.5$, $\lambda = 1$, $r=0.01$, $\alpha=0.01$ and the barrier of constant 
level $b$ is $10$. We plot the Laplace transform of ruin time and compare with 
the numerical results of collocation method. The results also illustrate the 
desireable performance of our proposed method.

\begin{figure}[h]
	\centering  %图片全局居中
	\subfigure[]{
		\label{Fig.sub.13}
		\includegraphics[width=0.45\textwidth]{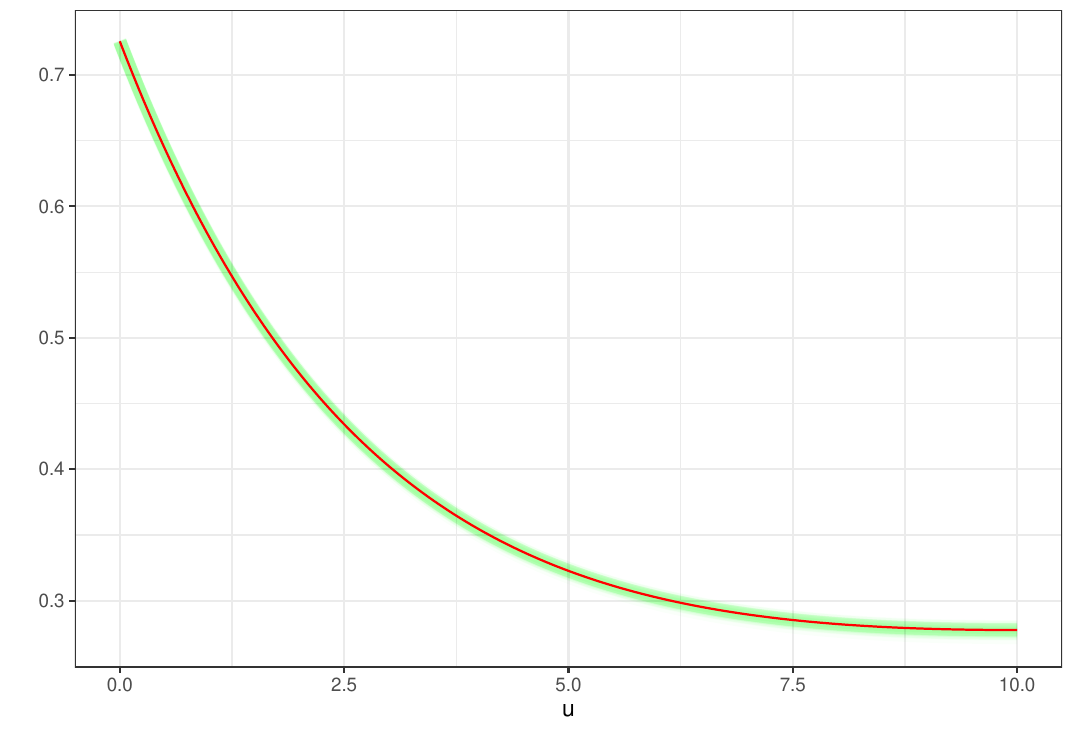}}
	\subfigure[]{
		\label{Fig.sub.14}
		\includegraphics[width=0.45\textwidth]{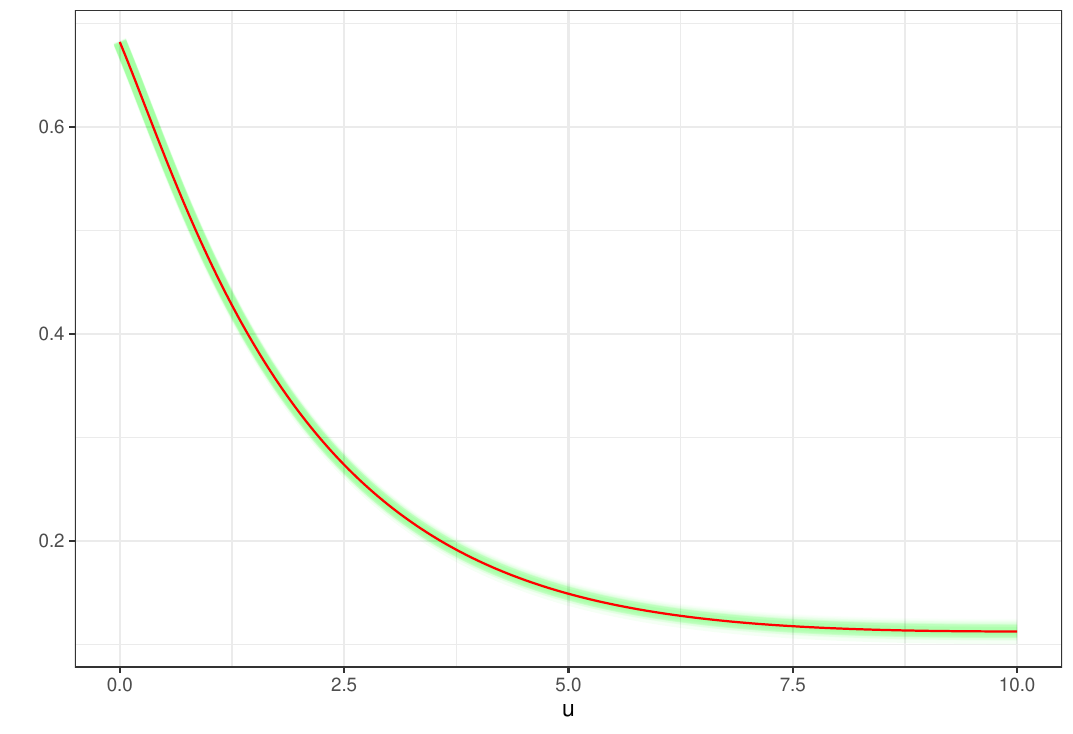}}
	\subfigure[]{
		\label{Fig.sub.15}
		\includegraphics[width=0.45\textwidth]{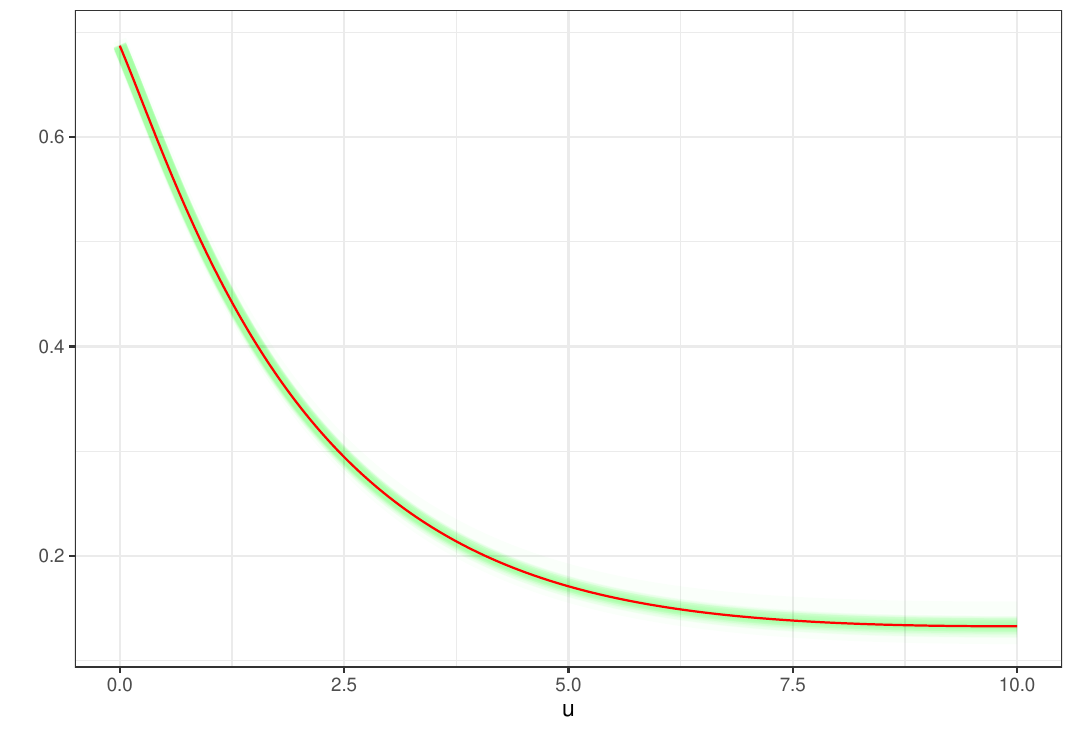}}
	\caption{Computing the Gerber-Shiu functions with a constant dividend 
	barrier. (a) exponential claim size density; (b) Erlang (2) claim size 
	density; (c) combination of exponentials}
	\label{Fig.5}
\end{figure}
	
\section{Concluding remarks}\label{sec:5}
In this paper, we have introduced the physics-informed neural networks to 
compute the Gerber-Shiu function. It is demonstrated that our algorithm is 
applicable when the Gerber-Shiu function contains interest and a constant 
dividend barrier. And this calculation method is also easy to implement in 
programming. 

We must note however that the proposed methods should not be viewed as 
replacements of classical numerical methods for computing the Gerber-Shiu 
functions, which, in many cases, meet the robustness and computational 
efficiency standards required in practice. For example, when dealing with the 
integral term on \eqref{2.5}, we use the Gaussian quadrature as an 
alternative. When $u$ becomes large, numerical integration will undoubtedly 
have a great impact 
on the calculation accuracy. Therefore, how to optimize the algorithm to 
improve calculation accuracy and efficiency is an issue that we will focus on 
in 
the future.

\appendix
	
%\begin{appendices}
	
\renewcommand{\thesection}{Appendix~\Alph{section}}
	
\renewcommand{\thesubsection}{\Alph{section}.\arabic{subsection}}
	
\renewcommand{\theequation}{\Alph{section}.\arabic{equation}}

\section{A representation of the discounted penalty function } 
\label{appendix:A}
\cite{lin2003classical} and \cite{yuen2007gerber} demonstrated that the 
Gerber-Shiu discounted penalty function with a barrier can be decomposed into 
two functions. The first function is the Gerber-Shiu discounted penalty 
function without a barrier. The second function, denoted as $h(u)$, is a 
nontrivial solution to the following homogeneous integro-differential equation:
\begin{equation}
		(r u+c) h^{\prime}(u)=(\lambda+\alpha) h(u)-\lambda \int_0^u h(u-y) 
		\mathrm{d} F(y), \quad u \geq 0 .
\end{equation}
Note that the value of $h(0)$ can be arbitrarily chosen due to the absence of 
constraint. Without loss of generality, we define the initial condition to be 
$h(0) = 1$.
	
Now consider the expected present value of the penalty at ruin in the absence 
of a dividend barrier, $\Phi_{\infty}(u)$. Similar to \eqref{2.5}, we have
\begin{equation}
		0=-(\alpha+\lambda) \Phi_{\infty}(u)+\Phi_{\infty}^{\prime}(u)(u 
		r+c)+\lambda 
		\int_0^u \Phi_{\infty}(u-y) \mathrm{d} F(y)+\lambda A(u) .
		\label{A.15}
\end{equation}
Thus, for $0 \leq u \leq b$, we have
\begin{equation}
		\Phi_{b}(u)=\Phi_{\infty}(u)-\frac{\Phi^{\prime}_{\infty}(b)}{h^{\prime}(b)}h(u)
\end{equation}
Hence, in order to determine $\Phi_{b}(u)$, we need to calculate both 
$\Phi_{\infty}(u)$ and $h(u)$, as well as their derivatives at $u = b$.
	
In fact, $\Phi_{\infty}(u)$ can be expressed as the following Volterra 
integral equation of the second kind, see \cite{cai2002expected}:
\begin{equation}
	\Phi_{\infty}(u)=\frac{c \Phi_{\infty}(0)}{c+ru}-\frac{\lambda}{c+r u} 
	\int_0^u A(t) \mathrm{d} t+\int_0^u 
	k(u, t) \Phi_{\infty}(t) \mathrm{d} t,
	\label{A.17}
\end{equation}
where
$$
k(u, t)=\frac{r+\alpha+\lambda \bar{F}(u-t)}{c+r u}.
$$ 
Similarly, $h(u)$ can be expressed as
$$
h(u)=\frac{c h(0)}{c+r u} +\int_0^u k(u, t) h(t) \mathrm{d} t.
$$
\section{The exact solution of $\Phi_{\infty}(0)$} \label{appendix:B}
The solution of \eqref{A.15} or \eqref{A.17} depends on the determination of 
the initial value, $\Phi_{\infty}(0)$. When $\alpha=0$, \cite{cai2002expected} 
give the 
specific experssion of $\Phi_{\infty}(0)$. In this part, we 
discuss the more general result when $\alpha>0$.

Define the following auxiliary function:
\begin{equation}
	Z_{\infty}(u)=\frac{\Phi_{\infty}(0)-\Phi_{\infty}(u)}{\Phi_{\infty}(0)}
	\label{B.1}
\end{equation}
Then $Z_{\infty}(0)=0$. Obviously, \eqref{B.1} implies that
\begin{equation}
	\Phi_{\infty}(u)=\Phi_{\infty}(0)-\Phi_{\infty}(0)Z_{\infty}(u).
	\label{B.2}
\end{equation}	
Subsititution of \eqref{B.2} into \eqref{A.17} results in
\begin{equation}
	(c+ru)Z_{\infty}(u)=(r+\alpha)\int_{0}^{u}Z_{\infty}(t) \mathrm{d}t
	-\alpha u+\frac{\lambda \mu_A}{\Phi_{\infty}(0)}A_1(t)
	-\lambda \mu F_1(u)+\lambda \mu F_1*Z_{\infty}(u),
	\label{B.3}
\end{equation}
where $A_1(u)=\frac{1}{\mu_A}\int_{0}^{u}A(t) \mathrm{d}t$, 
$\mu_A=\int_{0}^{\infty}A(t) \mathrm{d}t$, $F_1$ is the equilibrium 
distribution of $F$, given by
$$
F_1(u)=\frac{1}{\mu}\int_{0}^{u} 1-F(t) \mathrm{d}t,
$$
and
$F_1*Z_{\infty}(u)$	is the convolution of $F_1$ and $Z_{\infty}$.  
	
Define the Laplace transform of $Z_{\infty}(x)$, $F_1(x)$ and $A_1(x)$, namely,
$$
\tilde{z}_{\infty}(s)=\int_{0}^{\infty}e^{-sx}\mathrm{d}Z_{\infty}(x),
$$	
$$
\tilde{f_1}(s)=\int_{0}^{\infty}e^{-sx}\mathrm{d}F_{1}(x)
=\frac{1}{\mu}\int_{0}^{\infty}e^{-sx}\mathrm{d}F(x),
$$
$$
\tilde{a_1}(s)=\int_{0}^{\infty}e^{-sx}\mathrm{d}A_{1}(x)
=\frac{1}{\mu_A}\int_{0}^{\infty}e^{-sx}\mathrm{d}A(x).
$$

The Laplace transform of \eqref{B.3} with respect to $u$ results in the 
first-order differential equation for the function $\tilde{z}_{\infty}(s)$:
\begin{equation}
	c\tilde{z}_{\infty}(s)-r 
	\frac{\mathrm{d}}{\mathrm{d}s}\tilde{z}_{\infty}(s)=\frac{\alpha}{s}
	\tilde{z}_{\infty}(s)-\frac{\alpha}{s}+\frac{\lambda 
	\mu_A}{\Phi_{\infty}(0)}\tilde{a_1}(s)-\lambda\mu\tilde{f_1}(s)
+\lambda\mu\tilde{f_1}(s)\tilde{z}_{\infty}(s).
\label{B.4}
\end{equation}
Let 
$L(s)=c-\lambda\mu\tilde{f_1}(s)$,$M(s)=\frac{\lambda\mu_A}{\Phi_{\infty}(0)}
\tilde{a_1}(s)-\lambda\mu\tilde{f_1}(s)$, \eqref{B.4} is equivalent to 
$$
-r 
\frac{\mathrm{d}}{\mathrm{d}s}\tilde{z}_{\infty}(s)+\left(L(s)-\frac{\alpha}{s}
\right)\tilde{z}_{\infty}(s)=M(s)-\frac{\alpha}{s},
$$
When $r>0$, we note that
$$
\frac{\mathrm{d}}{\mathrm{d}s}\left(\tilde{z}_{\infty}(s)s^{\frac{\alpha}{r}}
\exp\left(-\frac{1}{r}\int_{0}^{s}L(t)\mathrm{d} t\right) \right)=-\frac{1}{r}
(M(s)-\frac{\alpha}{s})s^{\frac{\alpha}{r}}
\exp\left(-\frac{1}{r}\int_{0}^{s}L(t)\mathrm{d} t\right),
$$	
By monotone convergence we find that 
$$
\tilde{z}_{\infty}(s)s^{\frac{\alpha}{r}}=\int_{0}^{\infty}s^{\frac{\alpha}{r}}
e^{-sx}\mathrm{d}Z_{\infty}(x) \rightarrow 0, \quad (s \rightarrow \infty),
$$
which result in
\begin{equation}
	\tilde{z}_{\infty}(s)s^{\frac{\alpha}{r}}
	\exp\left(-\frac{1}{r}\int_{0}^{s}L(t)\mathrm{d} 
	t\right)=\frac{1}{r}\int_{s}^{\infty}
	(M(z)-\frac{\alpha}{z})z^{\frac{\alpha}{r}}
	\exp\left(-\frac{1}{r}\int_{0}^{z}L(t)\mathrm{d} t\right)\mathrm{d} z.
	\label{B.5}
\end{equation}	
When $\alpha>0$, substitute $s$ with zero in \eqref{B.5} to obtain
\begin{equation*}
\begin{aligned}
0=&\frac{\lambda \mu_A}{r 
		\Phi_{\infty}(0)}\int_{0}^{\infty}\tilde{a_1}(z)z^{\frac{\alpha}{r}}
	\exp\left({-\frac{1}{r}\int_{0}^{z}(c-\lambda\mu\tilde{f_1}(t))\mathrm{d} 
		t}\right)\mathrm{d} z \\
	&-\frac{1}{r}\int_{0}^{\infty}(\lambda\mu\tilde{f_1}(z)z+\alpha)z^{\frac{\alpha}{r}-1}
	\exp\left({-\frac{1}{r}\int_{0}^{z}(c-\lambda\mu\tilde{f_1}(t))\mathrm{d} 
		t}\right)\mathrm{d} z, 
\end{aligned} 
\end{equation*}
that is,
\begin{equation}
	\begin{aligned}
\Phi_{\infty}(0)&=\frac{\lambda \mu_A 
	\int_{0}^{\infty}\tilde{a_1}(z)z^{\frac{\alpha}{r}}
	\exp\left({-\frac{1}{r}\int_{0}^{z}(c-\lambda\mu\tilde{f_1}(t))\mathrm{d} 
		t}\right)\mathrm{d} 
		z}{\int_{0}^{\infty}(\lambda\mu\tilde{f_1}(z)z+\alpha)z^{\frac{\alpha}{r}-1}
	\exp\left({-\frac{1}{r}\int_{0}^{z}(c-\lambda\mu\tilde{f_1}(t))\mathrm{d} 
	t}\right)\mathrm{d} z}\\
&=\frac{\lambda 
\mu_A}{\kappa_{r,\alpha}}\int_{0}^{\infty}\tilde{a_1}(rv)v^{\frac{\alpha}{r}}
\exp\left(-cv+\int_{0}^{v}(\lambda\mu\tilde{f_1}(rs))\mathrm{d} 
s\right)\mathrm{d} v
\end{aligned}
\label{B.6}
\end{equation}
where \eqref{B.6} follows from the substitution $z=rv$, $t=rs$, and we define 
\begin{equation}
	\kappa_{r,\alpha}=c\int_{0}^{\infty}v^{\frac{\alpha}{r}}
	\exp\left(-cv+\int_{0}^{v}(\lambda\mu\tilde{f_1}(rs))\mathrm{d} 
	s\right)\mathrm{d} v.
	\label{B.7}
\end{equation}

\section*{Acknowledgement}
The author would like to thank the anonymous reviewers for their valuable 
suggestions which significantly improved the paper.

\bibliographystyle{apacite}

\bibliography{ref}
\addcontentsline{toc}{section}{References}

\clearpage
	
\end{document}